\documentclass[12pt, reqno]{amsart}
\usepackage{amsfonts}
\usepackage{amsmath}
\usepackage[applemac]{inputenc}
\usepackage[english]{babel}
\usepackage{amsthm}
\usepackage{graphicx}
\usepackage{color}
\usepackage{hyperref}

\usepackage{amssymb}
\usepackage{amsopn}
\usepackage{mathrsfs}
\usepackage{enumerate}

\usepackage{calrsfs}

\newcommand{\N}{\mathbb{N}}
\newcommand{\R}{\mathbb{R}}

\newcommand{\h}{\mathbb{H}}

\newcommand{\e}{\varepsilon}
\newcommand{\vp}{\varphi}
\newcommand{\pa}{\partial}
\newcommand{\p}{\partial}

\renewcommand{\P}{P_H}

\newcommand{\E}{E_{\mathrm{isop}}}  
\newcommand{\nablaH}{\nabla\!_H }

\newcommand{\res}
{\mathop{\hbox{\vrule height 7pt width .5pt depth 0pt \vrule
height .5pt width 6pt depth 0pt}}\nolimits}

\theoremstyle{plain}        
\theoremstyle{remark}        
\theoremstyle{definition}                      
\newtheorem{thm}{Theorem}[section]

\numberwithin{equation}{section}

\title[Quantitative isoperimetric inequalities
in $\h^n$]{Quantitative isoperimetric inequalities
in $\h^n$}
\author[V. Franceschi]{Valentina Franceschi }
\author[G.P. Leonardi]{Gian Paolo Leonardi}
\author[R. Monti]{Roberto Monti}
\date{\today}

\begin{document}
\maketitle

\linespread{1.2}

\begin{abstract}
In the Heisenberg group $\h^{n}$, $n\geq 1$,
we prove quantitative isoperimetric inequalities for Pansu's spheres, that  are
known to   be isoperimetric under various assumptions. The inequalities are
shown for suitably restricted classes of competing sets and the proof relies on
the construction of sub-calibrations.
\end{abstract}

\section{Introduction}

Quantitative isoperimetric inequalities in the Euclidean space and in
Riemannian manifolds have been an object of intensive studies in recent years.
The sharp quantitative isoperimetric
inequality in the Euclidean space $\R^{n}$ states that  there exists a constant
$C_n >0$ depending only on the dimension $n$, such that for any Borel set
$F\subset \R^{n}$ with $\mathcal L^n(F) = \mathcal L^n(B_1)$, the Lebesgue
measure of a unit ball $B_{1}$, one has the following estimate for the
difference of perimeters
\[
P(F) - P(B_{1}) \ge C_n \inf_{x\in \R^{n}}\mathcal L^n (F\Delta (x+B_{1}))^{2}.
\]
This inequality is established in its full generality in \cite{FMP}, and  
proved by different methods in \cite{FiMP, CL}. Several generalization   have
been recently obtained in Riemannian manifolds (with density), like the Gauss
space \cite{BarBraJul2014pre, CiaFusMagPra}, the $n$-dimensional sphere
\cite{BoeDuzFus2013pre}, and the $n$-dimensional hyperbolic space
\cite{BoeDuzSch2014pre}. A recurrent technique used in the proofs is based on
the regularity theory for perimeter quasiminimizers combined with a penalization
trick and a Fuglede-type argument, which essentially exploits the strict
positivity of the second variation of the area with respect to non-trivial
volume-preserving perturbations (see
\cite{CL, AceFusMor2013}). With  similar arguments, quantitative
stability results for global area-minimizing smooth hypersurfaces are obtained
in \cite{DePMag}, together with more specific results for a subfamily of
singular area-minimizing Lawson cones. In this case, due to the
presence of a singular point at the vertex of the cone, the proof of the sharp
quantitative stability follows   a different strategy, that is based on the
construction of suitable sub-calibrations (see also \cite{DePPao}).

On the other hand, in the context of subriemannian geometry, 
and in particular in Carnot groups, very few is known about the optimal constant
in the isoperimetric inequality (except for the fact that isoperimetric sets
exist and have at least some very weak regularity properties \cite{LR}). With
the only exception of the Grushin plane \cite{MonMor2004} (see also \cite{FM}),
isoperimetric sets have been only partially characterized in the subriemannian
Heisenberg group $\h^{n}$ (see below) and are not known at all in more general
Carnot groups.

The main result of this paper 
is Theorem \ref{THM1}, where we prove  the validity of a parameterized family of
quantitative isoperimetric inequalities in the subriemannian Heisenberg group
$\h^{n}$.

Before stating the result, 
we recall some basic definitions. The
$2n+1$-dimensional Heisenberg group  is the manifold
$\mathbb H^n=\mathbb C^n
\times\R$, $n\in\N$, endowed with the group product
\begin{equation*}
 \label{law!}
   (z,t)\ast (\zeta,\tau) = \big( z+\zeta, t+\tau+ 2 \, \mathrm{Im}\langle   z
, \bar \zeta\rangle \big),
\end{equation*}
where $t,\tau\in\R$, $z,\zeta\in \mathbb C^n$ and $\langle   z ,
\bar\zeta\rangle  =
z_1\bar\zeta_1+\ldots+z_n\bar\zeta_n$.
The bundle of horizontal left-invariant 
vector fields in $\mathbb H^n$ is
spanned by the vector fields
\begin{equation*}
 \label{XY}
  X_j = \frac{\p}{\p x_j}+2y_j\frac{\p}{\p t}, \quad
  Y_j = \frac{\p}{\p y_j}-2x_j\frac{\p}{\p t}
\end{equation*}
with $z_j = x_j+i y_j$ and $j=1,\ldots,n$.

The Haar measure of $\mathbb H^n$ is the 
Lebesgue measure $\mathcal L^{2n+1}$. The \emph{$H$-perimeter} of a $\mathcal L
^{2n+1}$-measurable set
$E\subset\mathbb H^n$ in an open set $A\subset\mathbb H^n$ is
\[
  \P(E,A) =
\sup\left\{ \int_E  \mathrm{div}\!_H  V
dzdt
:
V \in C^1_c(A ;\R^{2n}),\|V\|_\infty\leq 1 \right\},
\]
where 
  the horizontal divergence of the vector field
$V:A  \to \R^{2n}$    is
\[
\mathrm{div}\!_H   V=
\sum_{j=1}^n  X_jV_j+ Y_j V_{n+j}.
\]
We  use the notation $\mu_E(A ) =\P(E,A)$ and $\P(E) = \P(E,\mathbb H^n)$.
If $\P(E)<\infty$  then the open sets mapping $A\mapsto \mu_E(A)$  extends to a
Radon
measure $\mu_E$ on $\mathbb H^n$.
Moreover, there exists a $\mu_E$-measurable 
function $\nu_E:\mathbb H^n \to \R^{2n}$ such
that $|\nu_E| =1$ $\mu_E$-a.e.~and
the Gauss-Green integration by parts formula
\begin{equation}\label{Gauss-Green}
   \int_ {\mathbb H^n}  \langle V,
\nu_E\rangle   \, d\mu_E = - \int _{\mathbb H^n}  \mathrm{div}\!_H  
V\, dzdt
\end{equation}
holds for any $V \in C^1_c(\mathbb H^n; \R^{2n})$.
Here and hereafter, $\langle \cdot,\cdot\rangle$ denotes
the standard scalar product in $\R^{2n}$.

The isoperimetric problem in the Heisenberg group 
consists in minimizing $H$-perimeter of sets with 
a given fixed volume.
By homogeneity with respect to the dilations $(z,t)\mapsto(\lambda z, \lambda^2
,t)$ for $\lambda>0$, this is equivalent to prove 
existence, uniqueness, and classify the minimizers of 
the minimum problem
\begin{equation}\label{PIsop1}
\inf \left\{ \frac{\P(E,\mathbb H^n)}
{\mathcal L^{2n+1}(E)^{\frac{2n+1}{2n+2}}} \, : \,
\textrm{$E\subset\mathbb H^n$
measurable set with $0<\mathcal L ^{2n+1}(E)<\infty$} \right\}.
\end{equation}
A   set  realizing the
infimum is 
called \emph{isoperimetric set}. 
The existence of isoperimetric sets is established in \cite{LR}.

In 1983 P.~Pansu \cite{Pan1} conjectured that, up to left
translation and dilation,  the
 isoperimetric 
set is
\begin{equation} \label{E_isop}
\E 
   =\big\{
   (z,t)\in\h^n : |t|<\arccos(|z|)+|z|\sqrt{1-|z|^2},\,  |z|<1\big \}.
\end{equation}
The conjecture was made for  dimension $n=1$. 
The boundary of set $\E\subset\mathbb H^1$ can be obtained taking one geodesic  for the Carnot-Carath\'eodory metric joining the south pole
$(0,-\pi/2)\in\partial\E$
to the north pole $(0,\pi/2)\in\partial\E$ 
and letting it rotate around the $t$-axis.

In $\mathbb H^1$,  
Pansu's conjecture is proved 
assuming either the $C^2$ regularity of the minimizer \cite{RR}  or its 
convexity
\cite{MR}. In $\mathbb H^n$ with $n\geq 1$, the conjecture is proved 
assuming the axial symmetry of the
minimizer \cite{M3} or assuming a suitable cylindrical structure \cite{R}.  Some
observations on the problem can be found in \cite{M2} and \cite{LM}. See also the
book \cite{CDST} and the lecture notes \cite{M4}.

By refining the calibration argument 
of \cite{R} via a sub-calibration, we prove two  \emph{quantitative} versions of
the Heisenberg isoperimetric
inequality for competitors of $\E$ in half-cylinders.

For any $0\leq \e<1$ we define the
half-cylinder
\[
     C_\e = \big\{ (z,t) \in \mathbb H^n 
           : |z|<1 \text{ and }t>t_{\e} \big\},
\]
where $t_{\e} = \vp(1-\e)$  with 
$\vp(r)=\arccos(r)+r\sqrt{1-r^2}$.  
The proof provides an inequality with a variable structure, according to
whether 
$\e=0$ or $\e>0$. A similar 
construction could be used also in the Euclidean setting for Dido's problem
(i.e., for the relative isoperimetric problem in a half-space), and in this case
it would provide analogous quantitative estimates for the same classes of
competitors. Our main result is the following

\begin{thm}   \label{THM1}  
Let  $F\subset \mathbb
H^n$, $n\geq 1$, be any measurable set with  
$\mathcal L^{2n+1}(F) =\mathcal
L^{2n+1}(\E)$. 
\begin{itemize}
 \item [i)] If 
$F\Delta\E \subset\subset C_0$ then  
\begin{equation}
\label{TP}
\P(F) -\P(\E) \geq    \frac{n}{240\, \omega_{2n}^{2}}  \mathcal L^{2n+1}(
F\Delta\E)^3.
\end{equation}
 \item [ii)] If 
$F\Delta\E \subset\subset C_\varepsilon$ 
for $0<\varepsilon<1$, then  
\begin{equation}
\label{TP2}
\P(F) -\P(\E) \geq \frac{n\sqrt{\e}}
{16\,  \omega_{2n}}  \mathcal
L^{2n+1}(F\Delta\E)^2.
\end{equation}
\end{itemize}
Above,   $\omega_{2n}$ denotes the Lebesgue measure of the Euclidean unit ball
in
$\R^{2n}$.
\end{thm}

In \eqref{TP}, the asymmetry index 
$\mathcal L^{2n+1}(F\Delta\E)$ appears with the power $3$.
In \eqref{TP2}, the power is $2$ but there is a constant that vanishes with $\varepsilon$.
The quantitative isoperimetric inequality in $\mathbb R^n$ \cite{FMP} 
shows that the optimal power is $2$.

The sub-calibration is constructed in the following way.
The set $\E\cap C_\varepsilon$ can be foliated 
by a family of
hypersurfaces with constant $H$-mean curvature 
that decreases from $1$, the $H$-curvature  of 
$\partial \E$, to $0$, 
the curvature  of the surface $\{ t = 
t_\varepsilon\}$.
The velocity of the decrease depends on the parameter $\varepsilon$.
The horizontal unit normal to the leaves gives the sub-calibration.

The $H$-mean curvature is defined in the following way. Let
$\Sigma\subset\mathbb H^n$ be a hypersurface that is locally given by the zero
set of a function $ u\in C^1$ such that $|\nabla\!_H u|\neq 0$ on $\Sigma$,
where
\begin{equation} \label{acca}
 \nabla\!_H u = \big(X_1u,\ldots,X_nu, Y_1 u,\ldots, Y_n u\big)
\end{equation}
is  the horizontal gradient of $u$.
Then we define the $H$-mean curvature of $\Sigma$ at the point $(z,t) \in
\Sigma $ as
\begin{equation} \label{kappa}
   H_\Sigma (z,t) =  \frac{1}{2n} \mathrm{div}\!_H\Big(
   \frac{\nablaH u(z,t)}{|\nablaH u(z,t)|}\Big).
\end{equation}
The definition depends on a choice of sign. We shall work with orientable
embedded hypersurfaces and so we can  choose the positive sign, 
$H(z,t)\geq 0$.  Then, the boundary of $\E$ has constant $H$-mean curvature $1$.
For a set $E = \{(z,t)\in\mathbb H^n: u(z,t)>0\}$ the horizontal normal
$\nu_E$ in the Gauss-Green formula  \eqref{Gauss-Green} is given on $\partial
E$ by the vector
\[
 \nu_E = 
   \frac{\nablaH u }{|\nablaH u |}.
\]

The proof of Theorem \ref{THM1} relies on the construction described in the following result.

\begin{thm}  \label{THM3}
Let $0\leq \e<1$.
There exists a continuous function $u: C_\varepsilon\to\R$ with
level sets $\Sigma_s =\big\{ (z,t) \in C_\varepsilon : 
u(z,t) = s\big\}$,
$s\in \R$, such that:
\begin{itemize}
\item[i)] $u\in C^1(C_\varepsilon\cap
E_{\mathrm{isop}})
\cap C^1(C_\varepsilon\setminus \E)$ and $\nablaH u/|\nablaH u|$ is
continuously defined on $C_\e\setminus\{z=0\}$; 

\item[ii)] $\bigcup _{s>1} \Sigma_s =C_\varepsilon \cap \E$ and $\bigcup
_{s\leq 1} \Sigma_s =C_\varepsilon \setminus  \E$;

\item[iii)] $\Sigma_s$ is a 
   hypersurface of class $C^2$ with constant $H$-mean curvature
$H_{\Sigma_s}=1/s$ for 
$s>1$ and $H_{\Sigma_s} = 1$ for  $s\leq 1$;

\item[iv)] For any point $(z,\vp(|z|)-t) \in\Sigma_s$ with $s>1$ we have
\begin{equation}
 \label{H_s_2}
  1-H_{\Sigma_s}  (z,\vp(|z|)-t)\geq \frac{1}{20} t^2\quad 
\textrm{when }\varepsilon=0.
\end{equation}
and 
\begin{equation}
 \label{H_s_1}
  1-H_{\Sigma_s} (z,\vp(|z|)-t)\geq \frac{\sqrt \varepsilon}{4}t \quad
\textrm{when }0<\varepsilon<1,
\end{equation}
\end{itemize}
\end{thm}
The estimates \eqref{H_s_2} 
and \eqref{H_s_1} are the basis of the  
two inequalities \eqref{TP} and \eqref{TP2}, respectively.

\section{Proof of Theorem \ref{THM3}}

In $C_\varepsilon\setminus\E$, the leaves $\Sigma_s$ are
vertical translations of the top part of the boundary $\partial\E$. 
In $C_\varepsilon\cap \E$, the leaves
$\Sigma_s$ are constructed in the following way:
the surface $\partial \E$ is first dilated 
by a factor larger than $1$, and then it is translated
downwards in such a way that, after the two operations,
the sphere $\{ (z,t) \in \partial \E: t = t_\varepsilon\}$ 
with $t_\varepsilon = 
\vp(1-\varepsilon)$  remains fixed.

 The profile function of the set $\E$ is the
function $\vp:[0,1]\to\R$  
\begin{equation} \label{profile}
  \vp(r)=\arccos(r)+r\sqrt{1-r^2}   \qquad 0\leq r \le 1.
\end{equation}
Its first and second order 
derivatives  are 
\begin{equation}\label{D12}
\begin{split}
\varphi'(r)  = \frac{-2r^2}{\sqrt{1-r^2}}\quad \textrm{and}\quad 
\varphi''(r) = \frac{2r (r^2-2)}{(1-r^2)^{3/2}},\quad 0\leq r<1.
\end{split}
\end{equation}
Notice that $\varphi'''(0) = -4$. 
We also need the function $\psi:[0,1)\to\R$ 
\begin{equation}\label{funzg}
  \psi(r) = 2\vp(r) - r\vp'(r) 
       = 2\left(\frac{r}{\sqrt{1-r^{2}}} + \arccos(r)\right).
\end{equation}
Its derivative is 
\begin{equation}\label{derfunzg}
   \psi'(r) = \vp'(r) - r\vp''(r) 
   = \frac{2r^{2}}{(1-r^{2})^{  3/2}}.
\end{equation}

We start the construction of the function $u$. On the set
$C_\varepsilon\setminus\E$ we let 
\begin{equation} \label{pinco}
   u(z,t)=\varphi(|z|)-t +1,\quad (z,t)\in C_\varepsilon
      \setminus \E.
\end{equation}
Notice that $u(z,\vp(|z|)) = 1$ for all $|z|<1$. 
We define the function $u$ in the set  
\[
 D_{\e} =  C_\e\cap \E=
    \big\{ (z,t) \in \E: |z|<1-\e,\, 
      t_{\e}< t< \varphi(|z|)\big\}. 
\]
We use the short  
notation   $r = |z|$ and $r_{\e}=1-\e$.
Let   $F_\e:D_{\e}\times (1,\infty)\to\R$ be the function 
\[
   F_\varepsilon(z,t,s) = s^{2}\big(\vp(r/s) - \vp(r_{\e}/s)\big) + t_{\e} - t.
\]
We claim that for any point $(z,t) \in D_\e$ there exists 
a unique $s>1$ such
that $F_\e(z,t,s)=0$. In this case,
we can define the function  $u(z,t):D_{\e}\to \R$ letting 
\begin{equation}
\label{uimplicita}
   F_{\e}(z,t,s) = 0\quad \text{if and only if}\quad 
    s = u(z,t).
\end{equation}

We prove the claim. For any $(z,t)\in D_{\e}$ we have
\[
\lim_{s\to 1^{+}}F_{\e}(z,t,s) = \vp(r) - t >0.
\]
Moreover, with a second order Taylor expansion 
of $\vp$ based on \eqref{D12}
we see that
\[
\lim_{s\to \infty}F_{\e}(z,t,s) = t_{\e}-t <0.
\]
Since $s\mapsto F_\e(z,t,s)$ is continuous, this proves 
the existence of a
solution of $F_\e(z,t,s)=0$. 
By \eqref{funzg}, the derivative in $s$ of $F_\e$ is  
\begin{equation}
   \label{dsFeps}
     \pa_{s}F_{\e}(z,t,s) 
     = s \big(\psi(r/s) - \psi(r_{\e}/s)\big),
\end{equation}
and thus by \eqref{derfunzg} we deduce  that $\pa_{s}F_{\e}(z,t,s)  <0$. 
This proves the uniqueness.

We prove claim iii). Namely, we prove that for any point $(z,t)
\in \Sigma_s $ with $s>1$ and   $z\neq 0$,  the $H$-mean curvature of
$\Sigma_s$ at $(z,t)$ is
\begin{equation}\label{CL2}
 H_{\Sigma_s}(z,t) = -\frac{1}{2n} \mathrm{div}\!_H\Big(
   \frac{\nablaH u}{|\nablaH u|}\Big) = \frac 1 s.
\end{equation}
We are using definition \eqref{acca} with a minus sign in order to have a
positive curvature.
The claim when $s\leq 1$ is analogous because $\Sigma_s$ 
is a vertical translation of the top part of $\partial \E$.

By the implicit function theorem, the 
derivatives of $u$ can be computed from
the partial
derivatives of $F_\e$. 
Using $\pa_{x_{i}} r = {x_{i}}/{r}$ 
and $\pa_{y_{i}}r =  {y_{i}}/{r}$, with $i=1,\ldots,n$ and 
$z = (x_{1}+iy_{1},\dots,x_{n}+iy_{n})$, we find 
\begin{equation}
   \label{dxFeps}
   \pa_{x_{i}} F_{\e}(z,t,s) 
  = \frac{s x_{i}}{r}\vp'(r/s) 
\quad 
\text{and}
\quad
\pa_{y_{i}} F_{\e}(z,t,s) 
  = \frac{s y_{i}}{r}\vp'(r/s).
\end{equation}
Letting $s=u (z,t)$, thanks to \eqref{uimplicita}, 
\eqref{dsFeps}, \eqref{dxFeps},   and \eqref{D12}  we obtain
\begin{align}
\label{dxueps}
  \pa_{x_{i}} u(z,t) &
  = -\frac{\pa_{x_{i}} 
     F_{\e}(z,t,s)}{\pa_{s}F_{\e}(z,t,s)} 
  = \frac{2rx_{i}}{s\sqrt{s^{2}-r^{2}} 
      \big(\psi(r/s) - \psi(r_{\e}/s)\big)},
\\ 
\label{dyueps}
  \pa_{y_{i}} u(z,t) &= 
-\frac{\pa_{x_{i}} F_{\e}(z,t,s)}{\pa_{s}F_{\e}(z,t,s)} 
= 
\frac{2ry_{i}}{s\sqrt{s^{2}-r^{2}} \big(\psi(r/s) - \psi(r_{\e}/s)\big)},
\\ 
\label{dtueps}
  \pa_{t} u(z,t) &
  = -\frac{\pa_{t}F_{\e}(z,t,s)}{\pa_{s}F_{\e}(z,t,s)} 
  = \frac{1}{s\big(\psi(r/s) - \psi(r_{\e}/s)\big)},
\end{align}
and thus 
\begin{equation}
\label{dx-dt}
\pa_{x_{i}} u = 2x_{i} 
\frac{r}{\sqrt{s^{2}-r^{2}}}\pa_{t} u
\quad 
\text{and} 
\quad 
\pa_{y_{i}} u = 2y_{i} 
\frac{r}{\sqrt{s^{2}-r^{2}}}\pa_{t} u.
\end{equation}
It is then immediate to compute 
\[
\begin{split}
   X_{i}u & = \pa_{x_{i}} u + 2y_{i}\pa_{t}u 
          = \frac{2rx_{i} + 2y_{i}
  \sqrt{s^{2}-r^{2}}}{s\sqrt{s^{2}-r^{2}} 
    \big(\psi(r/s) - \psi(r_{\e}/s)\big)},
\\
   Y_{i}u & = \pa_{y_{i}} u - 2x_{i}\pa_{t}u 
  = \frac{2ry_{i} - 2x_{i}\sqrt{s^{2}-r^{2}}}
  {s\sqrt{s^{2}-r^{2}} 
 \big(\psi(r/s) - \psi(r_{\e}/s)\big)},
\end{split}
 \]
and  the squared length of the horizontal gradient 
of $u$ in $D_\e$  is 
\begin{align*}
   |\nablaH u |^{2} &
  = \sum_{i=1}^{n} (X_{i}u)^{2}+(Y_{i}u)^{2}
 \\ 
  &= \sum_{i=1}^{n} \frac{4r^{2}(x_{i}^{2}+y_{i}^{2}) 
 + 4(x_{i}^{2}+y_{i}^{2})(s^{2}-r^{2})}
   {s^{2}(s^{2}-r^{2}) 
  \big(\psi(r/s) - \psi(r_{\e}/s)\big)^{2}}
 \\ 
  &= \frac{4r^{2}}{(s^{2}-r^{2}) 
  \big(\psi(r/s) - \psi(r_{\e}/s)\big)^{2}}.
\end{align*}
Note that $|\nablaH u (z,t)| = 0$ if and only if $z=0$. 
So for any $(z,t) \in D_\e$ with $z\neq 0$ we have 
\begin{equation}
 \label{aieps}
a_{i}(z,t) = - \frac{X_{i}u}{|\nablaH u|} = \frac{rx_{i}+y_{i}\sqrt{s^{2}-r^{2}}}{rs} = 
  \frac{x_{i}}{s}+y_{i}\frac{\sqrt{s^{2}-r^{2}}}{rs}
\end{equation}
and
\begin{equation}\label{bieps}
b_{i}(z,t) = - \frac{Y_{i}u}{|\nablaH u|} = 
 \frac{ry_{i}-x_{i}\sqrt{s^{2}-r^{2}}}{rs} = 
  \frac{y_{i}}{s}-x_{i}\frac{\sqrt{s^{2}-r^{2}}}{rs}.
\end{equation}

If $(z,t) \in \E$ tends to $(\bar z,\bar t)\in 
\partial \E$ with $\bar t>0$ and $\bar z \neq 0$, then  $s = u(z,t)$ 
converges to $1$, and from \eqref{aieps} and \eqref{bieps} we see that 
\[
\lim_{(z,t)\to (\bar z,\bar t)} \frac{\nabla\!_H u(z,t)}
{|\nabla\!_H u(z,t)|} = - \Big( \bar x + \bar y 
\frac{\sqrt{1-|\bar z|^2}}{|\bar z|},\bar y - \bar x \frac{\sqrt{1-|\bar
z|^2}}{|\bar z|} \Big)=\frac{\nabla\!_H u(\bar z,\bar t)}
{|\nabla\!_H u(\bar z,\bar t)|},
\]
where the right hand side is computed using the definition
\eqref{pinco} of $u$. This ends the proof of claim i).

Claim ii) is clear. We prove claim iii).
The auxiliary function  $w(r,s) = {\sqrt{s^{2}-r^{2}}}/ {rs}$ satisfies 
\begin{equation}
\label{dxw}
\pa_{x_{i}} w 
   =\frac{x_{i}}{r}  \pa_{r}w + \pa_{x_{i}}u \, \pa_{s}w,\quad 
\pa_{y_{i}} w = \frac{y_{i}}{r} \pa_{r}w+ \pa_{y_{i}}u\, \pa_{s}w,\quad 
\pa_{s}w = \frac{r}{s^{2}\sqrt{s^{2}-r^{2}}}.
\end{equation}
By \eqref{aieps}, \eqref{bieps}, 
\eqref{dx-dt},  and \eqref{dxw}  we obtain
\begin{align*}
X_{i} a_{i} + Y_{i}b_{i} &
= \pa_{x_{i}}a_{i} 
+ 2y_{i} \pa_{t}a_{i} + \pa_{y_{i}}b_{i} 
- 2x_{i} \pa_{t}b_{i}
\\
&= \frac 1s - \frac{x_{i}}{s^{2}}\pa_{x_{i}}u 
+ y_{i}\Big(\frac{x_{i}}{r} \pa_{r}w  
 + \pa_{x_{i}}u\, \pa_{s}w  \Big) 
 + 2y_{i}\Big(-\frac{x_{i}}{s^{2}} \pa_{t}u 
 + y_{i}\pa_{s}w \, \pa_{t} u\Big)
\\ 
&
\quad 
+\frac 1s - \frac{y_{i}}{s^{2}}\pa_{y_{i}}u
 - x_{i}\Big(\frac{y_{i}}{r} \pa_{r}w 
 +\pa_{y_{i}}u\,  \pa_{s}w  \Big ) 
 -2x_{i}\Big(-\frac{y_{i}}{s^{2}} \pa_{t}u 
 - x_{i}\pa_{s}w \, \pa_{t} u\Big)
\\ 
&= \frac 2s - \frac{x_{i}\pa_{x_{i}}u 
+ y_{i}\pa_{y_{i}}u}{s^{2}} 
+ 2(x_{i}^{2}+y_{i}^{2})\pa_{s}w \, \pa_{t}u
\\ 
&= \frac 2s - \frac{x_{i}\pa_{x_{i}}u
 + y_{i}\pa_{y_{i}}u}{s^{2}} 
+ \frac{2r(x_{i}^{2}+y_{i}^{2})\pa_{t}u}
{s^{2}\sqrt{s^{2}-r^{2}}}= \frac 2s.
\end{align*}
Summing over $i=1,\ldots, n$ and dividing by $2n$, we obtain
\eqref{CL2}.


We prove claim iv). We fix a point $z$ with $|z|<1-\e$ 
and for $0\leq t<\vp(|z|)- t_\e$ we define the function
\begin{equation}\label{pippo}
    f_z (t) = u(z,\vp(|z|)-t) =  s = \frac{1}{H_{\Sigma_s}},
\end{equation}
where $s\geq 1$ is uniquely determined 
by $(z,\vp(|z|)-t) \in \Sigma_s$. 
The function $t\mapsto f_z(t)$ is increasing and $f_z(0)=1$

By \eqref{dtueps}, the function $f_z$ 
solves the differential equation
\[
f_z'(t) = -\pa_{t} u(z,\vp(|z|)-t) = \frac{1}{f_z(t) 
\big(\psi(r_{\e}/f_z(t)) -
\psi(r/f_z(t))\big)}
\]
for all $0< t <\vp(|z|)-t_\e$, and since, by 
\eqref{derfunzg}, $\psi$ is strictly increasing, $f_z$
solves the differential inequality
\[
   f_z'(t) \ge \frac{1}{f_z(t) \big(\psi(r_{\e}/f_z(t)) 
 - \pi\big)}.
\]
On the other hand, for any $s>1$ we have
\begin{equation} 
 \label{pop}
\begin{split}
  s\big(\psi(r_{\e}/s) - \pi\big) 
   & = s \int_{0}^{r_{\e}/s} \psi'(r)\, dr
  \\ 
 & = s \int_{0}^{r_{\e}/s} \frac{2r^{2}}{(1-r^{2})^{3/2}}\, dr
\\
&\le r_{\e} \int_{0}^{r_{\e}/s} \frac{2r}{(1-r^{2})^{3/2}}\, dr
\\ 
&= 2r_{\e}\Big((1-(r_{\e}/s)^{2})^{-1/2} - 1\Big)
\\ 
&\le \frac{2}{\sqrt{s-r_{\e}}}.
\end{split}
\end{equation}

In the case $\e=0$ we   have $r_{\e}=1$ and inequality \eqref{pop} 
reads 
\[
s\big(\psi(1/s) - \pi\big) \le \frac{2}{\sqrt{s-1}}. 
\]
Hence, the function $f_z$ satisfies the differential inequality 
\[
f_z'(t)  \geq \frac 12  \sqrt{f_z(t)-1}, \qquad t>0. 
\]
An integration with $f_z(0) =1$ gives $f_z(t) \geq 1+ {t^{2}}/{16}$, and
thus by the relation \eqref{pippo}  and by the bound $t<\pi/2$ we find
\[
 1-H_{\Sigma_s}(z,\vp(|z|)-t) = 1-\frac{1}{f_z(t)} \geq \frac{t^2}{16+t ^2} \geq
\frac{1}{20} t^2.
\]
This is claim \eqref{H_s_2}.

When  $0<\e<1$,  inequality \eqref{pop}  implies  
\[
s\big(\psi(r_{\e}/s) - \pi\big) \le \frac{2}{\sqrt{\e}},
\]
and thus $f_z'(t) \geq \sqrt{\e}/{2}$, that gives $f_z(t)\geq 1 + t
\sqrt{\e}/{2}$.
In this case, we find
\[
 1-H_{\Sigma_s} (z,\vp(|z|)-t)= 1 - \frac{1}{f_z(t)} \geq \frac{2\sqrt\e
t}{4+\pi}\geq \frac{\sqrt
\e}{4} t. 
\]
This is claim \eqref{H_s_1}. This finishes the proof of Theorem
\ref{THM3}.

\section{Proof of Theorem \ref{THM1}}
In this section, we prove the quantitative isoperimetric estimates 
\eqref{TP} and \eqref{TP2}.
\medskip 

Let $u:C_\varepsilon\to\R$, $0\leq\varepsilon<1$,
be the function given by Theorem \ref{THM3}
and let    $\Sigma_s=\{(z,t) \in C_\varepsilon :u(z,t)=s\}$ be the leaves of
the foliation, $s\in\R$. 
On $C_\varepsilon \setminus \{|z|=0\}$ we define the vector
field $X: C_\e\setminus \{|z|=0\}\to\R^{2n}$ by
\[
X=-\frac {\nabla\!_H u}{|\nabla\!_H u|}.
\]
Both $u$ and $X$ depend on $\e$. In particular, $X$ satisfies the following
properties:
\begin{itemize}
\item[i)] $|X|= 1$;
\item[ii)] for $(z,t)\in\partial \E \cap C_\varepsilon$ we have
$X(z,t)=-\nu_{\E}(z,t)$, the horizontal unit
  normal
to $\partial \E$.
\item[iii)] For any point $(z,t) 
\in \Sigma_s $, $s\in\R$, we have, 
\begin{equation} \label{INEQ}
 \frac{1}{2n}\mathrm{div}\!_H X(z,t) =H_{\Sigma_s}(z,t) \leq H_{\Sigma_0}=1.
\end{equation}
\end{itemize}

\medskip
We start the proof. 
Let $F\subset\mathbb H^n$ be a set with finite $H$-perimeter
such that $\mathcal L^{2n+1}(F)=\mathcal L^{2n+1}(\E)$ and
$F\Delta \E\subset \subset C_\varepsilon$. 
By Theorem 2.5 in \cite{FSSC}, we can without loss of generality 
assume that $\partial F$ is of class  $C^\infty$.
For $\delta>0$, let
$\E^\delta=\{(z,t)\in \E\ :\ |z|>\delta\}$.
By  \eqref{INEQ} and by the Gauss-Green formula 
\eqref{Gauss-Green}, 
we have
\begin{align*}
 \mathcal L^{2n+1}(\E^\delta\setminus F)
 &
 =\int_{\E^\delta\setminus F}1\,dzdt
\geq\int_{\E^\delta\setminus F}
\dfrac{\mathrm{div}\!_H X}{2n}\;dzdt
\\
\nonumber&
=\dfrac{1}{2n}\left\{\int_{
\partial F\cap \E^\delta}\langle X,\nu_{F}\rangle d\mu_{F}
-\int_{(\partial \E^\delta)\setminus F}\langle
X,\nu_{\E^\delta}\rangle d\mu_{\E^\delta}\right\}.
\end{align*}
Observe that
$\mu_{\E^\delta}=\mu_{\E}\res\{|z|>\delta\}+\mu_{\{|z|>\delta\}}\res \E$
and $\mu_{\{|z|>\delta\}}(\E)\leq C \delta^{2n-1}$.
Letting $\delta\to0^+$ and using the Cauchy-Schwarz inequality, we obtain
\begin{equation}
\label{eq:E-F}
\begin{split}
\mathcal L^{2n+1}(\E\setminus F)
 &
 \geq\dfrac{1}{2n}\left\{\int_{
\partial F\cap \E}\langle X,\nu_{F}\rangle d\mu_{F}
-\int_{(\partial \E)\setminus F}\langle
X,\nu_{\E}\rangle d\mu_{\E}\right\}
\\
&
\geq \frac{1}{2n}\left\{
\mu_{\E}(C_\varepsilon\setminus F)-\mu_F(\E)
\right\}
\\
&=\frac{1}{2n}\{\P(\E,C_\varepsilon\setminus F)-\P(F,\E)\}.
\end{split}
\end{equation}
By a similar computation we also have
\begin{align}
\label{eq:F-E}
\mathcal L^{2n+1}(F\setminus \E)&
=\int_{F\setminus \E}1\,dzdt
=\int_{F\setminus \E}
\dfrac{\mathrm{div}\!_H X}{2n}\;dzdt
\\
\nonumber&
=\dfrac{1}{2n}\left\{\int_{
\partial \E\cap F}\langle X,\nu_{\E}\rangle d\mu_{\E}
-\int_{(\partial F)\setminus \E}\langle
X,\nu_{F}\rangle d\mu_{F}\right\}
\\
\nonumber
&\leq \frac{1}{2n}\left\{
\mu_F(C_\varepsilon\setminus \E)-\mu_{\E}(F)
\right\}\\
&=\frac{1}{2n}\{\P(F,C_\varepsilon\setminus \E)-\P(\E,F)\}.
\end{align}

On the other hand,
\[\begin{split}
   \int_{\E\setminus F}\dfrac{\mathrm{div}\!_H X}{2n}\;dzdt&
  =\int_{\E\setminus F}
\Big(1+\Big(\dfrac{\mathrm{div}\!_H
X}{2n}-1\Big)\Big)\;dzdt\\
&=\mathcal L^{2n+1}(E\setminus F)-\int_{\E\setminus
F}\Big(1-\dfrac{\mathrm{div}\!_H X}{2n}\Big)\;dzdt\\
&=\mathcal L^{2n+1}(\E\setminus F)-\mathcal G(\E\setminus F),
\end{split}\]
where
\[
\mathcal G(\E\setminus F)=
\int_{\E\setminus F}\Big( 1-\dfrac{\mathrm{div}\!_H
X}{2n}\Big)\;dzdt.
\]
From  \eqref{eq:E-F} and  \eqref{eq:F-E}, we obtain   
\[
 \begin{split}
    \dfrac{1}{2n}\left\{\P(\E,C_\varepsilon\setminus F)-\P(F,\E)\right\}&
\leq\int_{\E\setminus F}\dfrac{\mathrm{div}\!_H X}{2n}\;dzdt
\\
&
=\mathcal L^{2n+1}(\E\setminus F)-\mathcal G(\E\setminus F)
\\
&=\mathcal L^{2n+1}(F\setminus \E)-\mathcal G(\E\setminus F)
\\
&\leq\dfrac{1}{2n}\left\{\P(F, C_\varepsilon\setminus
\E)-\P(\E,F)\right\}-\mathcal G(\E\setminus F),
\end{split}
\]
that is equivalent to 
\begin{equation} \label{F_E}
 \P(F)-\P(\E)\geq 2n \mathcal G(\E\setminus F).
\end{equation}

For any $z$ with $|z|<1-\e$, we define
the vertical sections
 $\E^z =\{t\in\R : (z,t) \in \E\}$ and
$F^z =
\{t\in\R:(z,t) \in F\}$. 
By Fubini-Tonelli theorem, we
have
\[
 \begin{split}
\mathcal G(\E\setminus F) &=\int_{\E\setminus F}
 \left(1-\dfrac{\mathrm{div}\!_H
X}{2n}\right)\,dzdt
\\
&
= \int _{\{|z|<1\}} \int_{\E^z\setminus F^z} \left(1-\dfrac{\mathrm{div}\!_H
X(z,t)}{2n}\right)dt \,dz.
\end{split}
\]
The function  $t\mapsto \mathrm{div}_H X(z,t)$ is increasing, and thus 
letting $m(z) = \mathcal L^1 (\E^z\setminus F^z)$, by monotonicity we obtain
\[
\begin{split}
 \mathcal G(\E\setminus F) & 
 \geq \int_{\{|z|<1\}}
\int_{\vp(|z|)-m(z)}^{\vp(|z|)}
\left(1-\dfrac{\mathrm{div}\!_H
X(z,t)}{2n}\right)dt \,dz
\\
&
= \int_{\{|z|<1\}}
\int_{0}^{m(z)}
\left(1-\frac{1}{f_z(t)} \right) dt \,dz,
\end{split}
\]
where $f_z(t) = u(z, \vp(|z|)-t)$ is the function introduced in \eqref{pippo}.

By \eqref{H_s_2}, when $\e=0$ the function $f_z$ satisfies the 
estimate $1-1/f_z(t) \geq t^2/20$, and by H\"older inequality we find
\begin{equation}\label{pix}
\begin{split}
\mathcal G(\E\setminus F) & 
    \geq \frac {1}{20} \int_{\{|z|<1\}}
    \int_{0}^{m(z)} t^2  dt \,dz 
\\
& =\frac{1}{60} \int_{\{|z|<1\}} m(z)^3    \,dz
 \\
&
\geq \frac{1}{60\omega_{2n}^2} \Big(\int_{\{|z|<1\}}
 m(z)    \,dz\Big)^3
\\
&
=\frac{1}{480\omega_{2n}^2} \mathcal L^{2n+1} (\E\Delta F)
^3.
\end{split}
\end{equation}
From \eqref{pix}  and \eqref{F_E} we obtain \eqref{TP}.

By \eqref{H_s_1}, when $0<\e<1$ the function $f_z$ satisfies the 
estimate $1-1/f_z(t) \geq \sqrt\e t/4$ and  we
find
\begin{equation}\label{pox}
\begin{split}
\mathcal G(\E\setminus F) & 
    \geq \frac {\sqrt\e}{4} \int_{\{|z|<1\}}
    \int_{0}^{m(z)} t  \, dt \,dz 
\\
& =\frac{\sqrt\e}{8} \int_{\{|z|<1\}} m(z)^2    \,dz
 \\
&
\geq \frac{\sqrt\e}{8\omega_{2n}} \Big(\int_{\{|z|<1\}}
 m(z)    \,dz\Big)^2
\\
&
=\frac{\sqrt\e}{32\omega_{2n}} \mathcal L^{2n+1} (\E\Delta F)
^2.
\end{split}
\end{equation}
From \eqref{pox} and  \eqref{F_E} we obtain  claim \eqref{TP2}.

\end{document}